\documentclass{amsart}

\theoremstyle{plain}  
\newtheorem{thm}{Theorem}[section]
\newtheorem{prop}[thm]{Proposition}

\theoremstyle{definition}

\newtheorem{prob}[thm]{Problem}

\theoremstyle{remark}
\newtheorem*{rem}{Remark}

\newcommand{\p}{\mathbf P}
\newcommand{\R}{\mathbf R}
\newcommand{\rad}{\mathrm{rad}\,}

\makeatletter 
\def\@setcopyright{} 
\def\serieslogo@{\@empty} 
\makeatother 

\begin{document}
\title[Duality of the second fundamental form]
{Dual varieties and the duality of \\the second fundamental form}
\author{Tohsuke Urabe}
\address{Department of Mathematics,
Tokyo Metropolitan University\\
Minami-Ohsawa 1-1,
Hachioji-shi, Tokyo, 192-03,
Japan}
\email{urabe@comp.metro-u.ac.jp}
\dedicatory{Dedicated to Professor Tzee-Char Kuo
on his sixtieth birthday}
\begin{abstract}
First, we consider a compact real-analytic irreducible subvariety
$M$ in a sphere and its dual variety $M^\vee$. 
We explain that two matrices of the second fundamental 
forms for both varieties $M$ and $M^\vee$
can be regarded as the inverse matrices of each other.
Also generalization in hyperbolic space is explained.
\end{abstract}
\maketitle

\section{Spherical case}
\label{sphere}
In this article I would like to explain main ideas 
in my recent results 
on duality of the second fundamental form.  
(Urabe\cite{{urabe;dual}}.)

Theory of dual varieties in the complex algebraic geometry is
very interesting.
(Griffiths and Harris~\cite{griffiths-harris;geo},
Kleiman~\cite{kleiman;enume},
Piene~\cite{piene;polar}, Urabe~\cite{urabe;polar},
Wallace~\cite{wallace;tangency}.)
Let $\p$ be a complex projective space of dimension $N$,
and $X\subset \p$ be a complex algebraic subvariety.
The set of all hyperplanes in $\p$ forms another projective space
$\p^\vee$ of dimension $N$, which is called the \emph{dual projective
space} of $\p$. 
The dual projective space $(\p^\vee)^\vee$ of $\p^\vee$ is 
identified with $\p$.
The closure in $\p^\vee$ of the set of tangent hyperplanes to $X$
is called the \emph{dual variety} of $X$, and is denoted by $X^\vee$.
We say that a hyperplane $H$ in $\p$ is tangent to $X$,
if we have a smooth point $p\in X$ such that $H$ contains
the embedded tangent space of $X$ at $p$.
It is known that the dual variety $X^\vee$ is again a complex
algebraic variety, and the dual variety $(X^\vee)^\vee$ of
$X^\vee$ coincides with $X$.

We would like to develop similar theory in the real-analytic
category.
(Obata~\cite{obata;gauss}.)

First, we fix the notations.
Let $N$ be a positive integer, and $L$ be a  vector
space of dimension $N+1$ over the real field $\R$.
A fixed positive-definite inner product on $L$ is denoted by
$(\ \ ,\ \ \ )$.
By $S=\{a\in L|(a,a)=1\}$ we denote the unit sphere in $L$.
The sphere $S$ has dimension $N$.

We consider a compact real-analytic irreducible subvariety $M$
in $S$.
We assume moreover that $M$ has only ordinary singularities as
singularities. 

We have to explain the phrase of ``ordinary singularity'' here.
Let $X\subset L$ be a real-analytic subset.
For every point $p\in X$ we can consider the germ $(X,p)$ of $X$
around $p$.
The germ $(X,p)$ is decomposed into into irreducible components.
By $\dim (X,p)$ we denote the dimension of the germ $(X,p)$.
The germ $(X,p)$ is said to be \emph{smooth}, if $(X,p)$ is
real-analytically isomorphic to $(\R^n,0)$ where
$n=\dim (X,p)$ and $0$ is a point of $\R^n$.
A point $p$ of $X$ is said to be smooth, if the germ $(X,p)$
is smooth.
We say that $X$ has an \emph{ordinary singularity} at $p\in X$,
if every irreducible component of $(X,p)$ is smooth.

Let $M_{smooth}\subset M$ be the set of smooth points $p\in M$
with $\dim (M,p)=\dim M$.
Under our assumption $M_{smooth}$ is dense in $M$.

For every point $p\in M_{smooth}$ the tangent space $T_p(M)$
of $M$ at $p$ is defined. 
Note in particular that $T_p(M)$ is not an affine subspace but
a vector subspace in $L$ passing through the origin.
The tangent space $T_p(M)$ has dimension equal to $\dim M$.
A point $q\in S$ is a normal vector of $M$ in $S$ at a point $p\in M$,
if $q$ is orthogonal to $p$ and $T_p(M)$.
We say that a point $q\in S$ is a normal vector of $M$ in $S$,
if $q$ is a normal vector of $M$ in $S$ at some point $p\in M$.
By $M^\vee$ we denote the closure in $S$ of the set
of normal vectors $a$ of $M$ in $S$ with $(a,a)=1$,
and we call $M^\vee\subset S$ the \emph{dual variety} of
$M\subset S$.
The dual variety $M^\vee$ has a lot of interesting properties.
However, $M^\vee$ is not a real-analytic subset in general.

\begin{prop}
\label{dense}
Under our assumption the dual variety $M^\vee$ contains
a dense smooth real-analytic subset whose connected components
have the same dimension.
\end{prop}

Let $X\subset S$ be a subset containing a dense smooth real-analytic
subset whose connected components
have the same dimension.
Obviously we can define the dual variety $X^\vee$ of $X$
by the essentially same definition as above.

\begin{thm}
\label{duality}
Under our assumption $(M^\vee)^\vee=M\cup\tau(M)$,
where $\tau:S\rightarrow S$ denotes the antipodal map $\tau(q)
=-q$.
\end{thm}

\begin{rem} 
Note that $M\cup\tau(M)$ is a
compact real-analytic subset only with ordinary singularities 
as singularities.
For any compact real-analytic subset in $L$ only with ordinary 
singularities as singularities, 
the irreducible decomposition is possible.
Therefore, $M$ is an irreducible component of $M\cup\tau(M)$, and
we can recover $M$ from $M\cup\tau(M)$.
\end{rem}

There exists an open dense smooth real-analytic subset $V$ of
$M^\vee$ such that for every point $q\in V$ there exists a point
$p\in M$ such that
\begin{enumerate}
\item $q$ is a normal vector of $M$ in $S$ at $p$, and
\item $p$ is a normal vector of $M^\vee$ in $S$ at $q$.
\end{enumerate}

Moreover, there exists an open dense smooth real-analytic subset $U$ of
$M$ such that for every point $p\in U$ there exists a point
$q\in V$ satisfying the same conditions 1 and 2 above.

Choose arbitrarily a pair $(q,p)$ of a smooth point $q\in M^\vee$ and
a smooth point $p\in M$ satisfying conditions 1 and 2, and fix it.

The second fundamental form of $M$ at $p$ in the normal direction $q$
$$\widetilde{II}:\: T_p(M)\times T_p(M)\longrightarrow \R$$
and the second fundamental form of $M^\vee$ at $q$ in the normal 
direction $p$
$$\widetilde{II}^\vee:\: T_q(M^\vee)\times T_q(M^\vee)
\longrightarrow \R$$
are defined.
We set
\begin{eqnarray*}
\rad\widetilde{II}&=&\{X\in T_p(M)|\mbox{For every }Y\in T_p(M),\:
\widetilde{II}(X,Y)=0\}\\
\rad\widetilde{II}^\vee&=&\{X\in T_q(M^\vee)|\mbox{For every }
Y\in T_q(M^\vee),\:
\widetilde{II}^\vee(X,Y)=0\}.
\end{eqnarray*}

\begin{thm}[Duality of the second fundamental form]
\label{second}
\ \newline\vspace*{-10pt}
\begin{enumerate}
\item $T_p(M)=\rad\widetilde{II}+(T_p(M)\cap T_q(M^\vee))$
(orthogonal direct sum)
\item $T_q(M^\vee)=\rad\widetilde{II}^\vee+(T_p(M)\cap T_q(M^\vee))$
(orthogonal direct sum)
\item $L=\R p+\rad\widetilde{II}+(T_p(M)\cap T_q(M^\vee))
+\rad\widetilde{II}^\vee+\R q$
(orthogonal direct sum)
\item Let $X_1, X_2,\ldots, X_r$ be an orthogonal normal basis of
$T_p(M)\cap T_q(M^\vee)$.
The matrix $(\widetilde{II}(X_i, X_j))$ is the inverse matrix of
$(\widetilde{II}^\vee(X_i, X_j))$.
\end{enumerate}
\end{thm}

Proposition~\ref{dense} is the most difficult part to show
in our theory.
Once we obtain Proposition~\ref{dense}, it is not difficult
to deduce Theorem~\ref{duality} applying analogous arguments
in complex projective algebraic geometry. 
Theorem~\ref{duality} and Theorem~\ref{second} can be shown 
through computation on Maurer-Cartan forms.
Theorem~\ref{second} seems to have a lot of applications
in theory of subvarieties in a sphere.

You can download my preprint~\cite{urabe;dual} containing verification
at
\begin{center}
\begin{tabular}{ll}
http://urabe-lab.math.metro-u.ac.jp/ & (Japanese)\\
http://urabe-lab.math.metro-u.ac.jp/DefaultE.html & (English).
\end{tabular}
\end{center}

\section{Hyperbolic case}
\label{hyperbolic}

We can consider similar situations in hyperbolic case.
(Obata~\cite{obata;gauss}.)

Let $L$ be a vector
space of dimension $N+1$ over the real field $\R$ as in 
Section~\ref{sphere}.
Now, we consider a non-degenerate inner product $(\ \ ,\ \ \ )$
on $L$ with signarutre $(N, 1)$.
By $S$ we denote one of the two connected components of the set 
$\{a\in L|(a, a)=-1\}$ in $L$.
The hyperbolic space $S$ has dimension $N$.

Also in this case we consider a compact real-analytic irreducible 
subvariety $M$ in $S$ only with ordinary singularities as
singularities.

Let $S^\vee=\{a\in L|(a, a)=1\}$.
Note that also $S^\vee$ is a smooth real-analytic connected variety
with dimension $N$.
However, $S\cap S^\vee=\emptyset$, 
and the metric on $S^\vee$ is not definite.
We can define the dual variety $M^\vee$ of $M$ 
as a subset of $S^\vee$
by the essentially same definition as above.
The dual variety $(M^\vee)^\vee$ of $M^\vee$ can be defined 
as a subset of $S$.

Proposition~\ref{dense} and Theorem~\ref{second} hold also in this
case without any modification. 
Theorem~\ref{duality} is replaced by the following brief theorem:

\begin{thm}
\label{duality2}
In hyperbolic case under our assumption $(M^\vee)^\vee=M$.
\end{thm}

\begin{prob}
Give genelarization of theory of dual varieties in 
$C^\infty$-category.
\end{prob}


\begin{thebibliography}{99}

\bibitem{griffiths-harris;geo}
Phillip Griffiths and Joseph Harris, \emph{Algebraic geometry and local
  differential geometry}, Ann. scient. \'{E}c. Norm. Sup. \mbox{$4^e$}
  s\'{e}rie \textbf{12} (1979), 355--432.

\bibitem{kleiman;enume}
Steven~L. Kleiman, \emph{The enumerative theory of singularities}, Real
and
  complex singularities (P.~Holm, ed.), Proceedings of the Nordic Summer
  School/NAVF, (Oslo, August 1976), Sijthoff \& Noordhoff, Alphen aan
den Rijn,
  The Netherlands, 1977, pp.~297--396.

\bibitem{obata;gauss}
Morio Obata, \emph{The \mbox{Gauss} map of immersions of
\mbox{Riemannian}
  manifolds in spaces of constant curvature}, J.~Differential Geometry
  \textbf{2} (1968), 217--223.

\bibitem{piene;polar}
Ragni Piene, \emph{Polar classes of singular varieties}, Ann. scient.
\'{E}c.
  Norm. Sup. \mbox{$4^e$} s\'{e}rie \textbf{11} (1978), 247--276.

\bibitem{urabe;polar}
Tohsuke Urabe, \emph{Duality of numerical characters of polar loci},
Publ.
  RIMS, Kyoto Univ. \textbf{17} (1981), 331--345.

\bibitem{urabe;dual}
Tohsuke Urabe, \emph{The Gauss map and the dual variety of real-analytic
submanufolds in a sphere or in hyperboloc space}, preprint (1995).

\bibitem{wallace;tangency}
Andrew~H. Wallace, \emph{Tangency and duality over arbitrary fields},
Proc.
  London Math. Soc. (3) \textbf{6} (1956), 321--342.

\end{thebibliography}
\end{document}